\documentclass[conference]{IEEEtran}
\IEEEoverridecommandlockouts
\usepackage{cite}
\usepackage{amsmath,amssymb,amsfonts}
\usepackage{algorithmic}
\usepackage{graphicx}
\usepackage{textcomp}
\usepackage{xcolor}
\def\BibTeX{{\rm B\kern-.05em{\sc i\kern-.025em b}\kern-.08em
    T\kern-.1667em\lower.7ex\hbox{E}\kern-.125emX}}

\usepackage[square,numbers]{natbib}

\usepackage[utf8]{inputenc}

\usepackage{hyperref}

\usepackage{subcaption}

\usepackage{dsfont}

\newcommand{\Eq}[1]{Equation~\eqref{#1}}

\newcommand{\Fig}[1]{Figure~\ref{#1}}

\newcommand{\dpart}[2]{\frac{\partial #1}{\partial #2}}

\newcommand{\dpartsq}[2]{\frac{\partial^2 #1}{\partial #2^2}}

\newcommand{\delpart}[2]{\frac{\delta #1}{\delta #2}}

\newcommand{\R}{\mathcal{R}}
\newcommand{\A}{\mathcal{A}}
\newcommand{\C}{\tilde{\mathcal{C}}}

\newcommand{\Lag}{\Lambda}

\newcommand{\ofx}{(x)}
\newcommand{\intL}{\int_{0}^{L}}
\newcommand{\lam}{\lambda}

\usepackage{xcolor}

\usepackage{amssymb}

\begin{document}

\title{Revisiting the optimal shape of cooling fins: \\ A one-dimensional analytical study using \\ optimality conditions
\thanks{A majority of this work took place while the first author was employed at the Technical University of Denmark.}
}

\author{\IEEEauthorblockN{Joe Alexandersen}
\IEEEauthorblockA{\textit{Department of Mechanical and Electrical Engineering} \\
\textit{University of Southern Denmark}\\
Odense, Denmark \\
joal@sdu.dk}
\and
\IEEEauthorblockN{Ole Sigmund}
\IEEEauthorblockA{\textit{Department of Mechanical Engineering} \\
\textit{Technical University of Denmark}\\
Kgs. Lyngby, Denmark}
}

\maketitle

\thispagestyle{plain}
\pagestyle{plain}

\begin{abstract}
This paper revisits the optimal shape problem of a single cooling fin using a one-dimensional heat conduction equation with convection boundary conditions.
Firstly, in contrast to previous works, we apply an approach using optimality conditions based on requiring stationarity of the Lagrangian functional of the optimisation problem. This yields an optimality condition basis for the commonly touted constant temperature gradient condition.
Secondly, we seek to minimise the root temperature for a prescribed thermal power, rather than maximising the heat transfer rate for a constant root temperature as previous works.
The optimal solution is shown to be fully equivalent for the two, which may seem obvious but to our knowledge has not been shown directly before.
Lastly, it is shown that optimal cooling fins have a Biot number of 1, exhibiting perfect balance between conductive and convective resistances.
\end{abstract}


\section{Introduction}

\citet{Duffin1959} presented a solution for the optimal shape of a single cooling fin to maximise the heat dissipation for a constant base temperature. His starting point was a statement from an earlier work by \citet{Schmidt1926}\footnote{The authors have not been able to find a copy of this work, hence base their citation purely on \citet{Duffin1959}.}: ``the temperature should be a linear function of the distance along the fin''. This is equivalent to stating that the heat flux should be constant along the distance of the fin. Based on rigorous mathematical analysis using calculus of variations, Duffin derived an optimal solution and proved the statement of Schmidt. Although Duffin provided a mathematical proof, we believe that a simpler physical explanation for the constant heat flux condition is still missing.

Our study is inspired by recent advances in topology optimisation for cooling design, reviewed by \citet{Dbouk2017} and \citet{Alexandersen2020}. The numerically obtained optimised designs generally have complex organic geometries and non-standard fin shapes, often resembling trees or corals \citep{Alexandersen2016,Alexandersen2018}. The recent analytical study by \citet{Yan2018} showed that contrary to popular belief, the organic tree-like structures are not optimal for pure conduction problems. Whether this also holds for the case of convection cooled heat sinks remains to be investigated.

In order to provide a starting point for subsequent studies, we revisit the problem using an approach which differentiates itself from the work of \citet{Duffin1959} in two ways: we attack the problem based on optimality conditions, requiring stationarity of the Lagrangian functional of the optimisation problem; and we treat a constant thermal power input at the root, rather than a constant root temperature.
Firstly, the optimality conditions present a simpler basis for the constant heat flux condition.
Secondly, although we tackle a different problem, it will be shown that the two are fully equivalent in the form of the final solution.

It should be noted that both Schmidt and Duffin made the customary assumption that the convective surface area was only dependent on the length and not the angle of the taper. \citet{Snider1987} argued against this assumption and postulated that a grooved fin surface must be the true optimum when the so-called ``length of arc'' is taken into account. However, the thermal performance of their solution is unbounded and prefers infinitely many infinitesimal grooves. This was formally proven by \citet{Marck2014}, who showed that there is no existence of solutions and that oscillating fins are preferred. Due to these issues of well-posedness, we rely on the same simple one-dimensional model as Duffin, with a convection surface area independent of the taper.

\section{One-dimensional cooling fin}

\begin{figure*}
    \begin{subfigure}{0.33\textwidth}
    \centering
    \includegraphics[width=\linewidth]{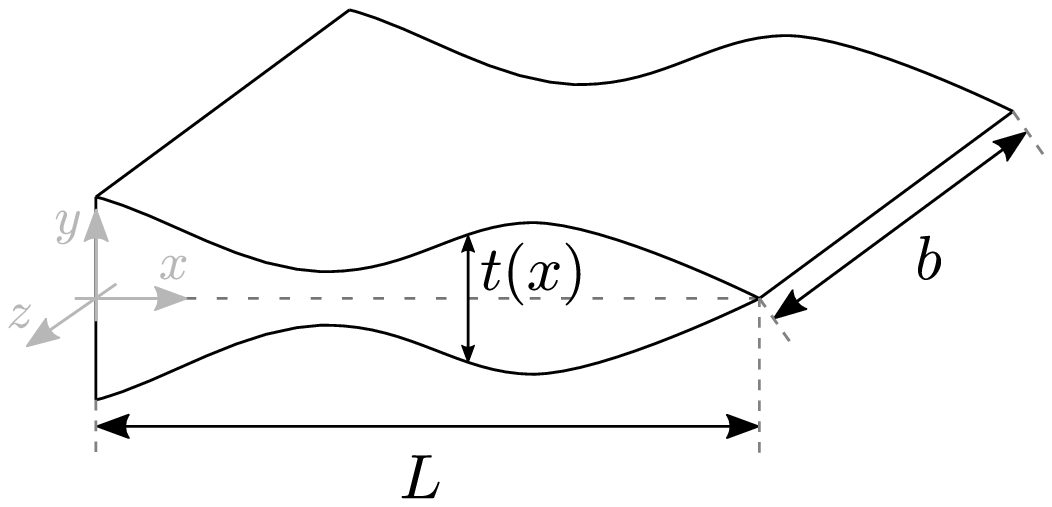}
    \caption{Dimensions}
    \label{fig:arbitraryFin-dims}
    \end{subfigure}
    \begin{subfigure}{0.33\textwidth}
    \centering
    \includegraphics[width=\linewidth]{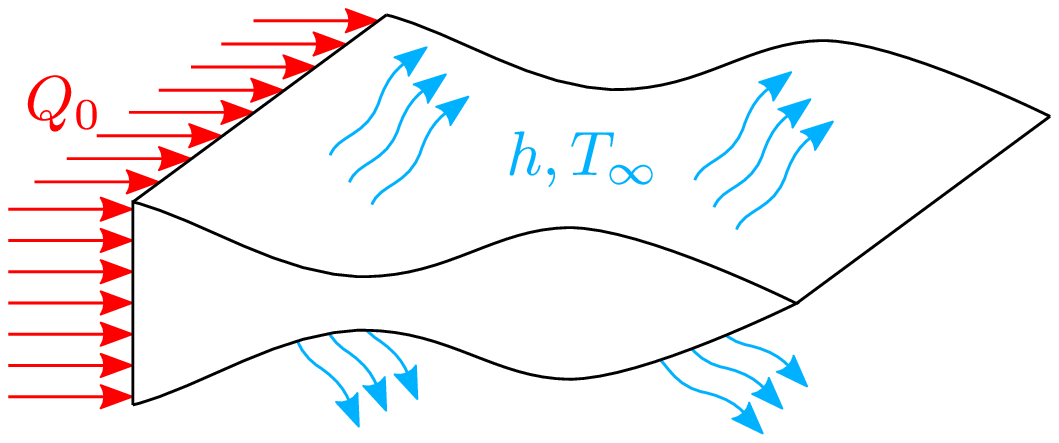}
    \caption{Boundary conditions}
    \label{fig:arbitraryFin-bcs}
    \end{subfigure}
    \begin{subfigure}{0.33\textwidth}
    \centering
    \includegraphics[width=\linewidth]{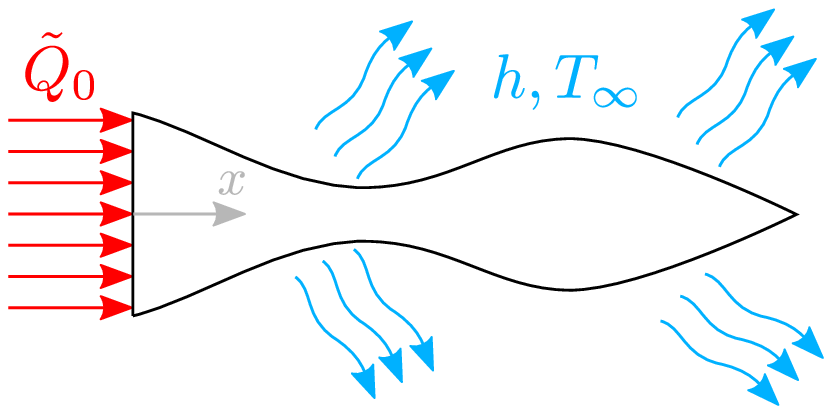}
    \caption{Simplification}
    \label{fig:arbitraryFin-plane}
    \end{subfigure}
    \caption{A fin of varying thickness subjected to a given thermal power, $Q_0$, distributed over the root surface and cooled by convection on the top and bottom surfaces characterised by the convection coefficient, $h$, and ambient temperature, $T_\infty$.}
    \label{fig:arbitraryFin}
\end{figure*}
A cooling fin of varying thickness is seen in \Fig{fig:arbitraryFin-dims}. The fin length is denoted $L$ and the fin width is denoted $b$. The fin profile in the $xy$-plane is characterised by the thickness as a function of the distance from the root, $t(x)$. The resulting profile is extruded in the negative $z$-direction, with the profile symmetric about the $x$-axis.
The thickness function, $t(x)$, is required to be strictly larger than zero in the domain, $t(x) > 0 \,\,\, \forall \,\,\, x \in [0;L[$, but is prescribed to end in a sharp tip, $t(L) = 0$. 

As shown in \Fig{fig:arbitraryFin-bcs}, the fin is subjected to a given constant thermal power input, $Q_0$, distributed over the root surface area. 
The fin is cooled on the top and bottom surfaces by a surrounding fluid through a convective flux expressed by Newton's law of cooling, $q_{c} = h ( T(x) - T_{\infty})$, where $h$ is the convection coefficient and $T_{\infty}$ is the ambient temperature.

It is assumed that no convection acts on the lateral ends (in the $z$-direction). Together with the assumption that the thickness at any point along the fin length is significantly smaller than the length and width, i.e. $t(x) << L$ and $t(x) << b$ for all $x \in [0;1]$, this allows for the assumption of a temperature field that only varies in the $x$-direction, i.e. $T(x,y,z) = T(x)$. This simplifies the problem to the one-dimensional problem seen in \Fig{fig:arbitraryFin-plane}.

By introducing the temperature difference,
\\ $\theta(x)~=~{T(x) - T_{\infty}}$, the governing equation with boundary conditions is given by:
\begin{subequations}\label{eq:govequ1}
\allowdisplaybreaks
\begin{align}
   - k \dpart{}{x} \left( t(x) \dpart{\theta}{x} \right) + 2 h \theta(x) &= 0 \,\,\, \text{for} \,\,\, 0 < x < L\label{eq:govequ1-a} \\
    -k\, t(0) \left.\dpart{\theta}{x}\right|_{x=0} &= \tilde{Q}_0  \label{eq:govequ1-b}\\
    -k\, t(L) \left.\dpart{\theta}{x}\right|_{x=L} &= 0 \label{eq:govequ1-c}
\end{align}
\end{subequations}
where $k$ is the thermal conductivity of the fin material and $\tilde{Q}_0 = \frac{Q_0}{b}$ is the thermal power input per unit width [W/m].
The above is derived from energy conservation for an infinitesimal slab under the assumption that the convection area is independent of the taper of the thickness profile.

\section{Optimisation problem}
By minimising the temperature at the root, we in turn maximise the heat transfer capabilities of the fin.
It is advantageous to define the objective functional as the temperature weighted by the applied thermal power per unit width, $\tilde{Q}_{0}$:
\begin{equation}
    \C = \theta(0) \cdot \tilde{Q}_{0}
\end{equation}
This is termed the thermal compliance, which is derived and related to the thermal resistance in Appendix \ref{app:thermalCompl}. This choice of objective functional causes the formulation to be self-adjoint as shown in Appendix \ref{app:selfadjoint}.

In order to ensure a non-trivial solution, an equality constraint is posed on the area of the profile (equivalent to the total volume of the fin):
\begin{equation}
    \A = A - \int_{0}^{L} t(x) dx = 0 \nonumber
\end{equation}
If there was no limit imposed on the material use, an infinite amount of material would be preferred.

Thus, the final optimisation problem becomes:
\begin{align} \label{eq:optprob}
    \underset{t \in \mathcal{D},L \in \mathds{R}_{+}}  {\mathrm{minimise}}  \,\,\,& \C \nonumber\\
    \mathrm{subject\;to} \,\,\,& \R = 0  \\
              \,\,\,& \A = 0 \nonumber 
\end{align}
where $\mathcal{D} = \left\lbrace t : \mathds{R} \rightarrow \mathds{R} \,\,|\,\, t(x) > 0 \, \, \forall \, x \in [0;L] \, , \, t(L) = 0 \right\rbrace$ is the admissible space of the thickness function and $\mathds{R}_{+} = \left\lbrace r \in \mathds{R} \,\,|\,\, r > 0 \right\rbrace$ is the admissible space of the length.

\section{Derivation of optimal solution}

Firstly, the Lagrangian functional is formed as:
\begin{equation}
    \Lag = \C - \intL w\ofx \R\ofx dx - \lam\A
\end{equation}
where $w$ and $\lam$ are Lagrange multipliers for the governing equation and area constraint, respectively, and $\R\ofx$ is the governing equations and boundary conditions recast as a residual, as detailed in Appendix \ref{app:deriveAnalytical}.

\subsection{Optimality conditions} \label{sec:optcond}

The first-order optimality conditions are given by requiring stationarity of the Lagrangian, $\Lag$, with respect to all variables:
\begin{equation}
    \delpart{\Lag}{\theta} = \delpart{\Lag}{w} = \dpart{\Lag}{\lam} = \delpart{\Lag}{t} = \dpart{\Lag}{L} = 0 \nonumber
\end{equation}
where $\dpart{\Lag}{v}$ is the partial derivative of $\Lag$ with respect to the variable $v$ and $\delpart{\Lag}{f}$ is the functional derivative of $\Lag$ with respect to the function $f$, as defined in Appendix \ref{app:unconstrDerivation}.

These stationarity conditions of the Lagrangian give rise to the following interesting optimality conditions on the solution:
\begin{subequations} \label{eq:optcrit}
\begin{align}
    w(x) &= \theta(x) & \forall x \in [0;L] \label{eq:optcrit-a}\\
    \dpart{\theta}{x} &= -\frac{\tilde{Q}_{0}}{h L^2} & \forall x \in [0;L]  \label{eq:optcrit-b}\\
    \dpart{t}{x} &= \frac{2h}{k} \left( x - L \right) & \forall x \in [0;L] \label{eq:optcrit-c}\\
    \theta(L) &= 0 \label{eq:optcrit-d}
\end{align}
\end{subequations}
These conditions lead to the following conclusions:
\begin{enumerate}
    \item \Eq{eq:optcrit-a} shows that the optimisation problem is self-adjoint, which simplifies the problem significantly.
    \item \Eq{eq:optcrit-b} shows that the temperature gradient is constant, which means that the heat flux is also constant. 
    \item \Eq{eq:optcrit-c} shows that the thickness gradient is linear and negative, which means that the thickness function has a convex quadratic taper.
    \item \Eq{eq:optcrit-d} shows that the temperature difference at the end is zero, which means that the optimal length is one that ensures all possible heat is lost. 
\end{enumerate}
These conclusions corroborate the arguments initially made by \citet{Schmidt1926} and subsequently proven mathematically by \citet{Duffin1959}.

\section{Optimal solution}

This section presents and discusses the solution, with the derivation detailed in Appendix \ref{app:deriveAnalytical}.

\subsection{Optimal length} 
The optimal length is given by:
\begin{equation} \label{eq:optsol-length}
    L = \sqrt[\leftroot{-1}\uproot{2}\scriptstyle 3]{\frac{3k A}{h}}
\end{equation}
which depends only on the allowable area, $A$, the conductivity, $k$, and the convection coefficient, $h$.
It can be seen that as the convection-to-conductivity ratio, $\frac{h}{k}$, increases, the optimal fin should be shorter. This result is corroborated by the well-known fact that as convection becomes stronger than conduction, a fin can be shorter because convection quickly removes the heat from it.

\subsection{Optimal thickness profile}

\begin{figure}
    \centering
    \includegraphics[width=\linewidth]{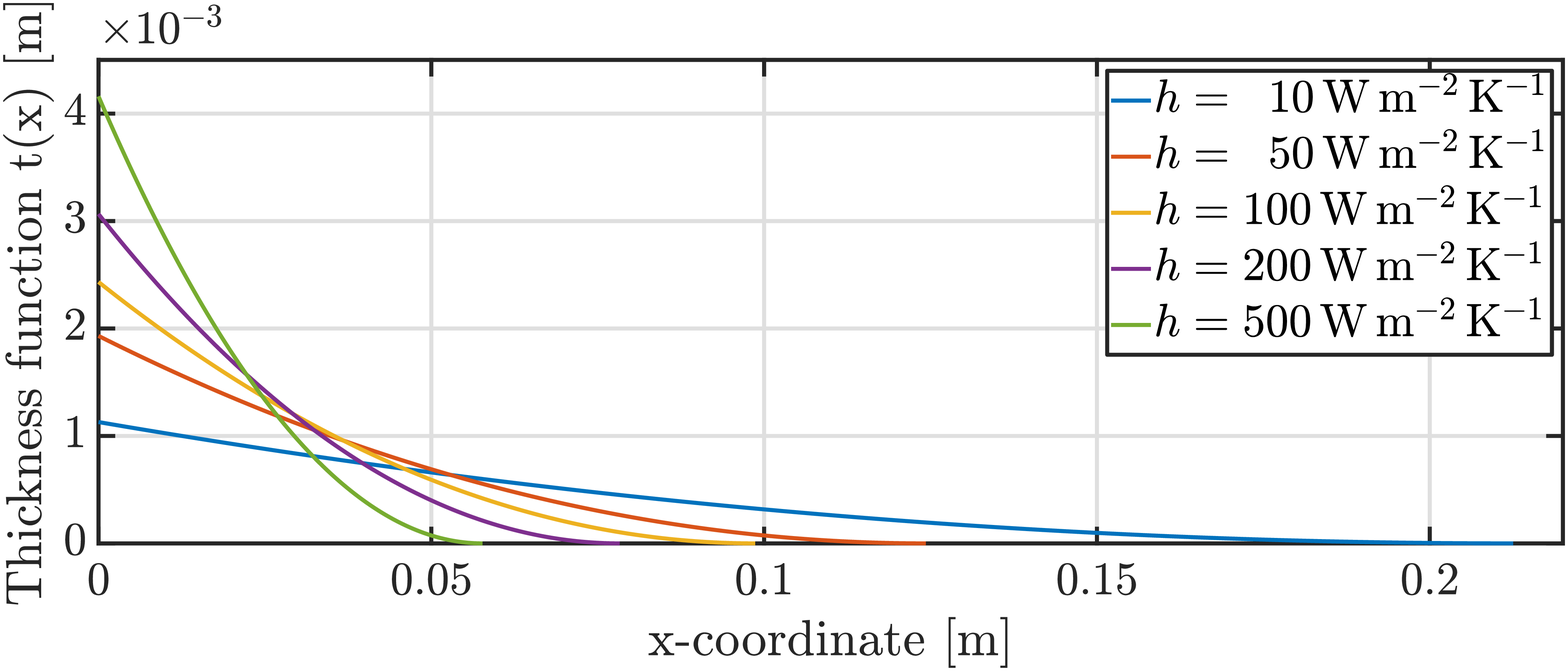}
    \caption{The optimal fin half-profiles for a range of convection coefficients for $k=200 \,\mathrm{W/(m\,K)}$, $A = 1.6\cdot10^{-4} \,\mathrm{m^{2}}$, and $\tilde{Q}_{0} = 20 \,\mathrm{W/m}$.}
    \label{fig:optimalThickness}
\end{figure}

The optimal thickness function is given by:
\begin{equation} \label{eq:optsol-thick}
    t\ofx = \frac{h}{k} \left( L - x \right)^2
\end{equation}
Figure \ref{fig:optimalThickness} shows the optimal thickness function for $k=200 \,\mathrm{W/(m\,K)}$, $A = 1.6\cdot10^{-4} \,\mathrm{m^{2}}$, and $\tilde{Q}_{0} = 20 \,\mathrm{W/m}$.
From \Eq{eq:optsol-thick} and Figure \ref{fig:optimalThickness}, it can be seen that the optimal fin has a convex and quadratic tapered profile. It is observed that: as the convection-to-conductivity ratio decreases, the optimal fin becomes thinner; as the convection-to-conductivity ratio increases, the optimal fin becomes wider.

\subsection{Optimal temperature profile}

\begin{figure}
    \centering
    \includegraphics[width=\linewidth]{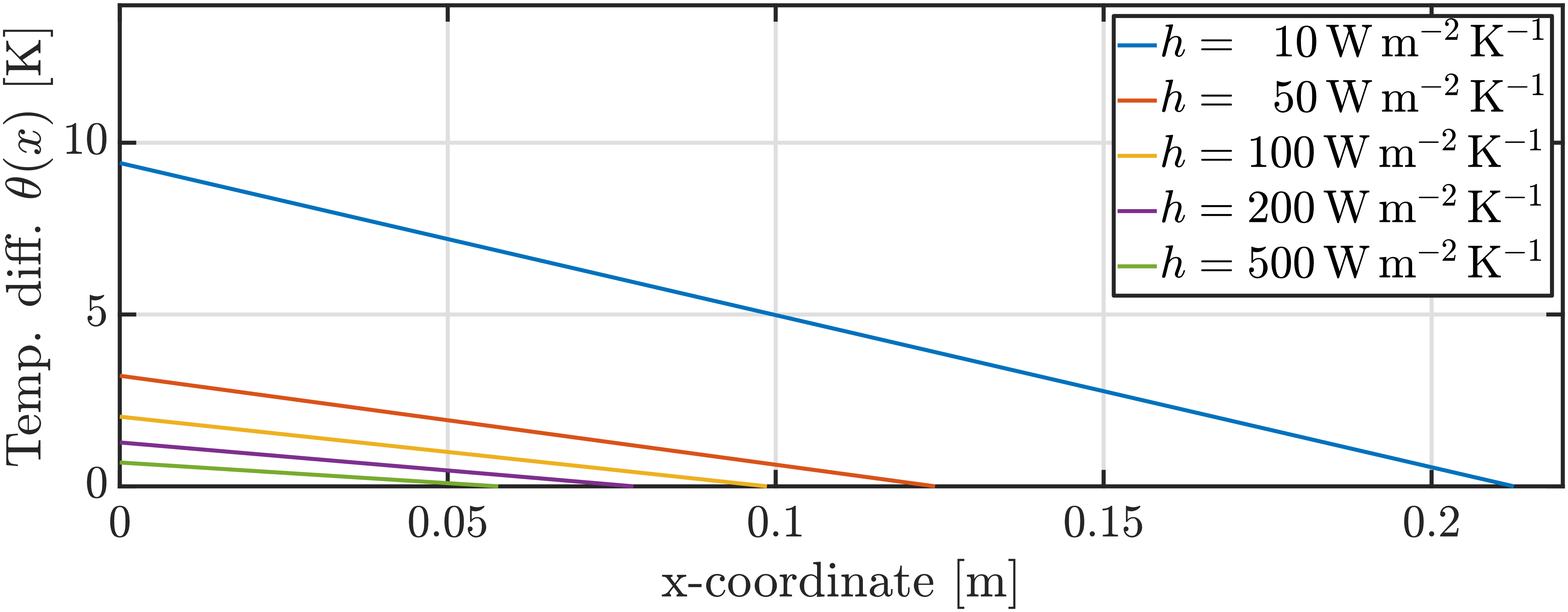}
    \caption{The optimal temperature profiles for a range of convection coefficients for $k=200 \,\mathrm{W/(m\,K)}$, $A = 1.6\cdot10^{-4} \,\mathrm{m^{2}}$, and $\tilde{Q}_{0} = 20 \,\mathrm{W/m}$.}
    \label{fig:optimalTemperature}
\end{figure}

The optimal temperature function is given by:
\begin{equation} \label{eq:optsol-temp}
    \theta\ofx = \frac{\tilde{Q}_{0}}{hL} \left( 1 - \frac{x}{L} \right)
\end{equation}
Figure \ref{fig:optimalTemperature} shows the optimal temperature function.
From \Eq{eq:optsol-temp} and Figure \ref{fig:optimalTemperature}, it is obvious that the temperature function is linear. This means that the optimal fin has a constant conductive heat flux through the fin.

It is evident from Figure \ref{fig:optimalTemperature} that the tip temperature difference is zero, as was also shown analytically by the optimality conditions in \Eq{eq:optcrit-d}. It makes sense physically that the tip temperature difference is zero, as this is the minimum obtainable value: the temperature cannot go below that of the ambient air. Thus, the optimal length of the fin is such that the temperature difference at the tip is zero. By having a shorter fin, the maximal amount of heat would not be dissipated by convection. By having a longer fin, the excess material would have no additional use.

\subsection{Equivalence to \citet{Duffin1959}}

By comparing the present solution, Equations~\eqref{eq:optsol-length} and \eqref{eq:optsol-thick}, to the solution derived by \citet[Equations (27) and (28)]{Duffin1959}, it can be seen that they are identical: \citet{Duffin1959} defines the thickness as $p(x) = \frac{1}{2}c(b-x)^2$ where the coefficient $c=\frac{2h}{k}$ and the length is given by $b=\sqrt[\leftroot{-1}\uproot{2}\scriptstyle 3]{\frac{6K}{c}}$ with $K=A$ being the constraint on the cross-sectional area.

It is well known that for a linear heat conduction problem, a given prescribed temperature boundary condition can be rewritten to one with a prescribed heat input yielding that temperature.
Applying a flux of $\tilde{Q}_{0} = hL = \sqrt[\leftroot{-1}\uproot{2}\scriptstyle 3]{ 3h^2 k A }$ ensures a root temperature difference of unity, $\theta_{0} = 1$, which is the problem solved by \citet{Duffin1959}. 

The two optimisation problems, maximising thermal energy throughput for a prescribed temperature and minimising the temperature for a prescribed thermal energy input, have now also been shown to be equivalent due to the same optimal solutions being found.

\subsection{Thermal compliance and resistance}

The optimal compliance is given by:
\begin{equation} \label{eq:optsol-compl}
    \C = \frac{\left.{\tilde{Q}_{0}}\right.^{2}}{hL} 
\end{equation}
As derived in Appendix \ref{app:thermalCompl}, the total thermal resistance per unit width of the optimal fin is given by:
\begin{equation} \label{eq:optsol-thres}
     \tilde{R}_\mathrm{fin} = \frac{\C}{\left.{\tilde{Q}_{0}}\right.^{2}} = \frac{1}{hL}
\end{equation}
The total thermal resistance is comprised of the conductive and convective contributions, $\tilde{R}_\mathrm{fin} = \tilde{R}_\mathrm{cond} + \tilde{R}_\mathrm{conv}$. The convective resistance per unit width is defined as $\tilde{R}_\mathrm{conv} = \frac{1}{2hL}$. For the unconstrained case of \Eq{eq:optsol-thres}, it is therefore evident that the conductive and convective resistances must be the same. If defined as the ratio of the conductive to convective resistances, this implies that the Biot number will always be 1 for the optimal fin. This means that the conductive and convective resistances are perfectly balanced.

\section{Further study}

This short study provides a basis for future studies into the optimal design of cooling fins and heat sinks. Revisiting the problem based on optimality conditions provides a clear reasoning for the optimal fin shape.

The conclusion that the Biot number should always be 1 for an optimal fin, should be investigated in practise for topology-optimised designs. Calculations of the Biot number by \citet{Alexandersen2016} for entire heat sinks have previously shown this not to be the case in general, possibly due to geometrical constraints and finite distances to neighbouring fin features. Whether the conclusion even extends to the case of the true (taper-dependent) convection surface area must also be investigated.

Furthermore, studies into whether truly optimal cooling fins and heat sinks exhibit branching behaviour for convective problems are needed. The ill-posedness and preference to infinitely many infinitesimal grooves \citep{Snider1987,Marck2014} may indicate that branching is preferable. But formal studies should be carried out to extend the knowledge from pure conduction \citep{Yan2018} to convection cooled heat sinks.

\section*{Acknowledgment}

The authors thank Anton Evgrafov for feedback on an early draft of this work.

\bibliographystyle{unsrtnat}      
\bibliography{fin1D.bib}   

\appendix

\subsection{Definition of thermal compliance}\label{app:thermalCompl}

The thermal compliance is defined as the integral over the fin root surface of the temperature multiplied by the heat flux:
\begin{equation}
    \mathcal{C} = \int_{A_{rs}} \theta_{0} \phi_0 \,dS
\end{equation}
where the applied heat flux is given by:
\begin{equation}
    \phi_0 = \frac{{Q}_{0}}{t_{0}b}
\end{equation}
Since both the heat flux and temperature are assumed constant at the entire root surface, with a surface area of $A_{rs} = t_{0}b$, the thermal compliance becomes:
\begin{equation}
    \mathcal{C} = t_{0}b \theta_{0} \phi_{0}  = {Q}_{0} \theta_{0}
\end{equation}

The thermal resistance of the fin is defined as the ratio between the difference between root and ambient temperatures to the thermal energy input, which can be related to the thermal compliance:
\begin{equation}
    R_\mathrm{fin} = \frac{\theta(0)}{{Q}_{0}} = \frac{\mathcal{C}}{{{Q}_{0}}^2}
\end{equation}
Since ${Q}_{0}$ is a constant, this means that minimising the thermal compliance is directly equivalent to minimising the thermal resistance of the fin.

\subsection{Derivation of optimal solution} \label{app:deriveAnalytical}

Firstly, the Lagrangian functional is formed as:
\begin{equation}
    \Lag = \C - \intL w\ofx \R\ofx dx - \lam\A
\end{equation}
where $w$ and $\lam$ are Lagrange multipliers for the governing equation and area constraint, respectively. The Lagrange multiplier for the residual is a spatially-varying field and is henceforth called the adjoint field. The governing equations and boundary conditions have been recast as a residual:
\begin{align}
   \R\ofx = & - k \dpart{}{x} \left( t(x) \dpart{\theta}{x} \right) + 2 h \theta(x) \nonumber \\
            & + \delta(x-0) \left( - k\, t(x) \dpart{\theta}{x} - \tilde{Q}_0 \right) \nonumber \\
            & + \delta(x-L) \left( - k\, t(x) \dpart{\theta}{x} \right) = 0
\end{align}
where $\delta(x-a)$ is the Dirac delta function active at $x=a$.

\subsubsection{Optimality conditions} \label{app:unconstrDerivation}

The first-order optimality conditions are given by requiring stationarity of the Lagrangian, $\Lag$, with respect to all variables:
\begin{equation}
    \delpart{\Lag}{\theta} = \delpart{\Lag}{w} = \dpart{\Lag}{\lam} = \delpart{\Lag}{t} = \dpart{\Lag}{L} = 0
\end{equation}
where $\dpart{\Lag}{v}$ is the partial derivative of $\Lag$ with respect to the variable $v$ and $\delpart{\Lag}{f}$ is the functional derivative of $\Lag$ with respect to the function $f$.
The functional derivative is defined from the following relation:
\begin{equation}
    \delta \Lag[f] = \varepsilon \intL \frac{\delta \Lag}{\delta f} \eta\, dx
\end{equation}
with the variation given by:
\begin{equation}
    \delta \Lag[f] = \Lag[f+\varepsilon\eta] - \Lag[f]
\end{equation}
where $\eta$ is a given function perturbation and $\varepsilon$ is the perturbation size.

The individual stationarity conditions yields the following set of equations:
\begin{subequations} \label{eq:statcrit}
\begin{align}
    \delpart{\Lag}{\theta} &= k \dpart{}{x} \left( t(x) \dpart{w}{x} \right) - 2 h w(x) -  \delta(x-L) \left( k t(x)\dpart{w}{x} \right)  \nonumber \\ &+ \delta(x-0)\left( k t(x) \dpart{w}{x} + \tilde{Q}_{0} \right) = 0 \,\,\,\, \forall\,\, x\in[0;L] \label{eq:statcrit-a} \\
    \delpart{\Lag}{w} &= k \dpart{}{x} \left( t(x) \dpart{\theta}{x} \right) - 2 h \theta(x) -  \delta(x-L) \left( k t(x)\dpart{w}{x} \right) \nonumber \\ &-\delta(x-0) \left( k\, t(x) \dpart{\theta}{x} + \tilde{Q}_0 \right) = -\R\ofx = 0 \,\,\,\, \forall\,\, x\in[0;L] \label{eq:statcrit-b} \\
    \dpart{\Lag}{\lam} & = A - \intL t(x) dx= \A = 0 \label{eq:statcrit-c} \\
    \delpart{\Lag}{t} &= - k\dpart{w}{x}\dpart{\theta}{x} + \lambda = 0 \,\,\,\, \forall\,\, x\in[0;L] \label{eq:statcrit-d} \\
    \dpart{\Lag}{L} &= \left.\left( \dpart{w}{x} k t(x) \dpart{\theta}{x} \right)\right|_{x=L} + 2 w(L) h \theta(L) - \lam t(L) = 0  \label{eq:statcrit-e}
\end{align}
\end{subequations}

\subsubsection{Self-adjointness} \label{app:selfadjoint}
Stationarity wrt. the temperature field, \Eq{eq:statcrit-a}, is fulfilled if $w\ofx$ satisfies the following:
\begin{subequations}\label{eq:adjointproblem}
\begin{align}
  - k \dpart{}{x}\left( t(x) \dpart{w}{x} \right) + 2 h w\ofx & = 0 & \mathrm{for} \,\,\, 0 < x < L \label{eq:adjointproblem-a} \\
    -k t(0) \dpart{w}{x} &= \tilde{Q}_{0} \label{eq:adjointproblem-b} \\ 
    -k t(L) \dpart{w}{x} &= 0 \label{eq:adjointproblem-c} 
\end{align}
\end{subequations}
This defines the adjoint problem and it can easily be seen that it is identical to \Eq{eq:govequ1} and thus:
\begin{equation}\label{eq:selfadjoint}
    w(x) = \theta(x)
\end{equation}

If the governing equation was posed in terms of the absolute temperature, $T\ofx$, the problem would not be self-adjoint because a $2 h T_\infty$ term would be missing on the right-hand side of Equation \ref{eq:adjointproblem-a}. However, there would be the simple relation that $w(x) = T\ofx-T_\infty$.

\subsubsection{Constant temperature gradient}
Stationarity wrt. the thickness field, \Eq{eq:statcrit-d}, is satisfied if:
\begin{equation}
    \dpart{w}{x}k\dpart{\theta}{x} = \lambda  \,\,\, \forall \,\,\, x \in [0;L]\nonumber \nonumber
\end{equation}
Introducing the fact that the problem is self-adjoint, \Eq{eq:selfadjoint}, this yields:
\begin{equation}
    \left( \dpart{\theta}{x} \right)^2 = \frac{\lam}{k}\nonumber
\end{equation}
Taking the square root of the above yields two solutions for the gradient, one positive and negative. Due to physical considerations, the thermal gradient must be negative and thus:
\begin{equation}\label{eq:constgrad}
    \dpart{\theta}{x} = - \sqrt{\frac{\lam}{k}} \nonumber
\end{equation}
This condition shows that the optimal fin is one with a constant temperature gradient. Therefore, it is known that the temperature field must be linear:
\begin{equation}\label{eq:lineartemp1}
    \theta\ofx = \theta_{0} - \sqrt{\frac{\lam}{k}}x
\end{equation}
where $\theta_{0}$ is the temperature at the root to be found next.

\subsubsection{Temperature at the root}
Performing the definite integral on stationarity wrt. the adjoint field, \Eq{eq:statcrit-b}, yields:
\begin{equation}
    -\intL \R\ofx dx =  \tilde{Q}_{0} - \intL 2 h \theta\ofx dx  = \nonumber
\end{equation}
This equations states that the heat energy sent into the root is all lost through convection along the fin.
Inserting the linear temperature field, \Eq{eq:lineartemp1}, into the above allows for explicit integration giving:
\begin{equation}
    2 h \left[ \theta_{0}x - \frac{1}{2}\sqrt{\frac{\lam}{k}}x^{2} \right]^{L}_{0} = \tilde{Q}_{0} \nonumber
\end{equation}
which when evaluated gives the expression for the root temperature:
\begin{equation}
    \theta_{0} = \frac{\tilde{Q}_{0}}{2hL} + \frac{L}{2}\sqrt{\frac{\lam}{k}} \nonumber
\end{equation}
The temperature field is then updated to:
\begin{equation}\label{eq:lineartemp2}
    \theta\ofx = \frac{\tilde{Q}_{0}}{2hL} + \sqrt{\frac{\lam}{k}} \left( \frac{L}{2} - x \right)
\end{equation}

\subsubsection{Tip thickness and tip temperature}
Since it has been defined that $t(L) = 0$, the first and third term of \Eq{eq:statcrit-e} are zero. This leaves the second term:
\begin{equation}
    2w(L)h \theta(L) = 0 \nonumber
\end{equation}
By introducing the self-adjointness, \Eq{eq:selfadjoint}, it is easily seen that:
\begin{equation} \label{eq:tiptemp}
    \theta(L) = 0
\end{equation}
which means that the optimal length is that giving a temperature difference of zero at the tip. This makes sense physically, as this is the minimum obtainable value. By having a shorter fin, the maximal amount of heat would not be dissipated by convection. By having a longer fin, the excess material would have no additional use.

\subsubsection{Lagrange multiplier for area constraint}
By inserting the temperature function, \Eq{eq:lineartemp2}, into the zero tip temperature condition, \Eq{eq:tiptemp}, yields the following expression for the Lagrange multiplier:
\begin{equation} 
    \lambda = \frac{k {\tilde{Q}_{0}}^2}{h^{2} L^{4}} \nonumber
\end{equation}
This expression is then used to update the temperature function to:
\begin{equation}\label{eq:lineartemp2}
    \theta\ofx = \frac{\tilde{Q}_{0}}{hL} \left( 1 - \frac{x}{L} \right)
\end{equation}
where it can be seen that the root temperature is simply given by:
\begin{equation}
    \theta_{0} = \frac{\tilde{Q}_{0}}{hL} \nonumber
\end{equation}

\subsubsection{Thickness function}
Using the product rule for differentiation on the diffusive term of the governing equation, \Eq{eq:govequ1-a}, gives:
\begin{equation}
   - k \left( \dpartsq{\theta}{x} t(x) + \dpart{\theta}{x}\dpart{t}{x} \right) + 2h \theta(x) = 0 \nonumber
\end{equation}
Inserting the temperature expression, \Eq{eq:lineartemp2}, gives:
\begin{equation}
   -k \frac{\tilde{Q}_{0}}{h L^{2}} \dpart{t}{x} - \frac{2\tilde{Q}_{0}}{L} \left( 1 - \frac{x}{L} \right) = 0 
\end{equation}
because the second-order derivative of the linear temperature function is zero.
From this, the gradient of the thickness can be isolated as:
\begin{equation}
    \dpart{t}{x} = \frac{2h}{k} \left( x - L \right)
\end{equation}
which can be indefinitely integrated to find the following expression for the thickness function:
\begin{equation}
    t\ofx = t_{0} + \frac{2h}{k} \left( \frac{x^{2}}{2} - xL \right)
\end{equation}
where $t_{0}$ is the introduced constant, here the thickness at the root. It can be found using the zero tip thickness condition, $t(L)=0$, giving:
\begin{equation}
    t_{0} = \frac{h L^{2}}{k}
\end{equation}
Thus, the final thickness function expression becomes:
\begin{equation} \label{eq:thickfunc}
    t\ofx = \frac{h}{k} \left( L - x \right)^2
\end{equation}

\subsubsection{Optimal fin length}
Finally, the optimal fin length can be found from stationarity wrt. the Lagrange multiplier for the area constraint. Inserting the thickness function, \Eq{eq:thickfunc}, into the area constraint gives:
\begin{equation}
    \intL \frac{h}{k} \left( L - x \right)^2 dx = A \nonumber
\end{equation}
which when integrated with bounds inserted yields the optimal length:
\begin{equation}
    L = \sqrt[\leftroot{-1}\uproot{2}\scriptstyle 3]{\frac{3k A}{h}} \nonumber
\end{equation}

\subsubsection{Final optimal solution} \label{sec:unconstrSolution}
The final optimal solution is:
\begin{subequations} \label{eq:finalOptSol}
\begin{align}
    L &= \sqrt[\leftroot{-1}\uproot{2}\scriptstyle 3]{\frac{3k A}{h}}  \label{eq:finalOptSol-a}\\
    t\ofx &= \frac{h}{k} \left( L - x \right)^2  \label{eq:finalOptSol-b}\\
    \theta\ofx &= \frac{\tilde{Q}_{0}}{hL} \left( 1 - \frac{x}{L} \right)  \label{eq:finalOptSol-c}\\
    \C &= \frac{\left.{\tilde{Q}_{0}}\right.^{2}}{hL} \label{eq:finalOptSol-d}
\end{align}
\end{subequations}

\end{document}